\documentclass[a4paper,11pt,twoside]{article}

\usepackage[T2A]{fontenc}
\usepackage[cp1251]{inputenc}
\usepackage[english,russian]{babel}
\usepackage{amsmath,amsfonts,amssymb,mathrsfs,amscd}
\usepackage{verbatim}
\usepackage{eucal}
\usepackage{graphicx}
\usepackage{enumerate}
\usepackage[pdftex,unicode,colorlinks,linkcolor=blue,
citecolor=red,bookmarksopen,pdfhighlight=/N]{hyperref}
\newtheorem{theorem}{\bf Теорема}%%[section]
\newtheorem{propos}{\bf Предложение}%%[section]
\newtheorem{proof}{\rm Доказательство}
\newtheorem{remark}{Замечание}%%[section]
\newcommand{\RR}{\mathbb{R}}
\newcommand{\CC}{\mathbb{C}}
\newcommand{\DD}{\mathbb{D}}
\newcommand{\const}{{\rm const}}

\renewcommand{\leq}{\leqslant}
\renewcommand{\geq}{\geqslant}

\DeclareMathOperator{\dist}{dist}
\DeclareMathOperator{\clos}{clos}
\DeclareMathOperator{\Int}{int}
\DeclareMathOperator{\Hol}{Hol}
\DeclareMathOperator{\CHol}{CHol}

\DeclareMathOperator{\Zero}{Zero}
\DeclareMathOperator{\sbh}{sbh}
\DeclareMathOperator{\shg}{shg}
\DeclareMathOperator{\trc}{trc}
\DeclareMathOperator{\bd}{bd}
\DeclareMathOperator{\ext}{ext}
\title{О множествах неединственности для пространств голоморфных функций}
\author{Б.\,Н.~Хабибуллин, Ф.\,Б.~Хабибуллин}

\begin{document}

\maketitle

\begin{abstract}
В статье  усовершенствуются и уточняются некоторые наши результаты из предшествующей нашей недавней статьи из журнала <<Известия вузов. Математика>> 2015 г. за счет последних наших результатов 2016 г. об оценках снизу субгармонических функций логарифмом модуля голоморфной ненулевой функции.
\end{abstract}

\section*{Введение}

Как обычно,  $\RR$ и  $\CC$
 --- множества соответственно всех {\it  вещественных\/} и  {\it комплексных чисел\/} или их естественные геометрические интерпретации; кроме того, $\DD:=\{z\in \CC \colon |z|<1\}$ --- {\it единичный круг\/} на  комплексной плоскости $\CC$. Используются  определения и понятия из \cite{Levin56}--\cite{Kh12}, но при необходимости мы их повторяем.
  Пусть $D$ --- область в $\CC$.
 Каждой  не более чем счетной последовательности точек $\Lambda =\{\lambda_k\}_{k\geq 1}\subset D$ без точек сгущения в $D$  сопоставляем {\it считающую меру} $n_\Lambda$, а именно:   
$	n_\Lambda (S):=\sum_{\lambda_k\in S} 1$ 
--- число точек из  $\Lambda$, попавших в $S\subset D$. Заметим, что  среди точек $\lambda_k$ могут быть и повторяющиеся. %%Объединение $\Lambda\cup \Lambda_0$ двух таких  последовательностей $\Lambda, \Lambda_0\subset D$ полностью %%определяется считающей мерой $n_{\Lambda\cup \Lambda_0}:=n_{\Lambda}+n_{\Lambda_0}$. 
По определению функция $n_\Lambda(\lambda):=n_\Lambda \bigl(\{\lambda\}\bigr)$ --- дивизор последовательности $\Lambda$, т.\,е. число повторений точки $\lambda\in \CC$ в последовательности $\Lambda$. Так, $\lambda \in \Lambda$, если $n_{\Lambda}(\lambda)>0$. 

Векторное пространство всех голоморфных в $D$ функций обозначаем как $\Hol (D)$. Если не оговорено противное, пространство $\Hol (D)$ наделяем топологией равномерной сходимости на компактах из $D$.
 Ненулевой функции $f\in \Hol (D)$ соответствует {\it последовательность нулей\/} $\Zero_f$, перенумерованная с учетом кратности.

Для компакта $C\subset \CC$  через $\CHol[C]$ обозначаем векторное пространство над полем $\CC$ непрерывных на $C$ комплекснозначных функций, одновременно голоморфных во {\it внутренности\/}  $\Int C$, если она не пуста, с естественной $\sup$-нормой.
 
 Последовательность точек $\Lambda\subset D$ называется {\it подпоследовательностью нулей\/} для подмножества $H\subset \Hol (D)$, если найдется ненулевая функция $f\in H$, для которой $\Lambda \subset \Zero_f$ в том смысле, что $n_\Lambda (\lambda)\leq  n_{\Zero_f}(\lambda)$ для всех $\lambda \in D$. Если  $H$ замкнуто относительно вычитания, например, векторное подпространство над $\RR$, то подпоследовательность нулей для $H$ называют и {\it последовательностью,\/} или {\it множеством, неединственности\/} для $H$. 
 
Выпуклый конус всех субгармонических фун\-кций в области $D\subset \CC$ обозначаем через $\sbh (D)$. Субгармоническую функцию, тождественно равную $-\infty$ на $D$, обозначаем $\boldsymbol{-\infty}$. Для $s\in \sbh (D)$ меру Рисса функции $s$ чаще всего будем обозначать как $\nu_s$, и наоборот, субгармоническую функцию $s$ в $D$ с мерой Рисса $\nu$ часто записываем в виде в виде $s:=s_\nu$. Борелевскую положительную меру  (конечную на компактах из $D$), или меру Радона $\nu$ \cite[Appendix A]{Rans}, называем {\it подмерой для подмножества\/} $S\subset \sbh (D)$, если найдется функция $s\in S$, которая $\neq \boldsymbol{-\infty}$, с мерой Рисса $\nu_s\geq \nu$ на $D$. Иначе говоря, $\nu$ --- подмера для $S$, если для некоторой (любой) субгармонической функции $s_\nu$ с мерой Рисса $\nu$ найдется функция $v\in \sbh (D)$, не равная $\boldsymbol{-\infty}$, для которой $s:=s_{\nu}+v\in S$. Возможность варьирования слов <<некоторый>> и <<любой>> в последнем предложении обеспечена <<нечуствительностью>> неравенств к перекидыванию гармонических слагаемых от одного субгармонического слагаемого к другому. 

Для (весовой) функции $M\colon D\to [-\infty,+\infty]$ со значениями в  расширенной вещественной оси $[-\infty,+\infty]:=\{-\infty\}\cup\RR \cup\{+\infty\}$ с естественным отношением порядка определим {весовой класс} субгармонических функций
\begin{equation*}
\sbh(D;M]:=\{s\in \sbh (D)\colon s\leq M +\const \; \text{\it  на $D$} \},
\end{equation*}
где здесь и далее $\const$ --- какая-либо постоянная, а <<{\it $s\leq M+\const$ на\/ $D$\,}>> означает выполнение поточечных неравенств $s(z)\leq M(z)+\const$ во всех точках $z\in D$.
Аналогично, определим весовое пространство голоморфных функций
\begin{equation*}
\Hol(D;\exp M]:=\{f\in \Hol (D)\colon |f|\leq \const \cdot \exp M \;  \text{\it  на $D$} \}.
\end{equation*}
 
 В   разделе \ref{sec2} рассматривается следующая задача. Пусть $N$ и $M$ --- две весовые функции в области $D\subset \CC$
и $\Lambda$ --- последовательность точек в $D$. При каких простых соотношениях между $N$ и $M$ некоторая подмера $\mu \geq n_\Lambda$  для $\sbh(D;M]$ определяет подпоследовательность нулей  $\Lambda$ для пространства $\Hol(D;\exp N]$ или, возможно,чуть большего пространства? 
Более или менее удовлетворительное решение этой задачи позволяет свести исследование подпоследовательностей нулей к гибкому аппарату субгармонических функций, к тому же в классах $\sbh(D;M]$, отличных от $\sbh(D;N]$. Также в важном подразделе \ref{gd} тот же вопрос отдельно исследуется для весовых пространств функций на всей плоскости $\CC$, определяемых положительно однородными при показателе $\rho>0$ весовыми функциями.

\section{Последовательности неединственности}\label{sec2}

Пусть $S$ --- подмножество расширенной комплексной плоскости $\CC_{\infty}:=\CC\cup \{\infty\}$.
Для подмножества $S\subset \CC_{\infty}$ через $\clos S$ и  $\bd S$ обозначаем соответственно 
замыкание и границу $S$ в $\CC_{\infty}$;
на евклидовом пространстве $\dist (\cdot , \cdot)$ --- евклидово расстояние между двумя объектами (точками, подмножествами) в евклидовом пространстве (в нашем случае $\RR$ или  $\CC$). 
Пусть $D$ --- область в $\CC$. Пусть $d\colon D\to (0,1]$ --- непрерывная функция, удовлетворяющая условию
\begin{equation}\label{ddis}
0<d(z)<\dist(z, \bd D), \quad z\in D.
\end{equation}  
Каждой  функции $N$ сопоставляем ее усреднение по кругам с центром $z$ радиуса $0<r<\dist(z,\bd D)$, обозначаемое как 
\begin{equation}
B  (z,r;N):=\frac{1}{\pi r^2}\int_0^{2\pi}\int_0^r N(te^{i\theta}) t\, d t \,d \theta.
\end{equation}

%%Пишем $D\Subset \CC$, если область $D$ {\it предкомпактна\/} в $\CC$, т.е. просто ограничена.
Весовой функции $N\colon D\to [-\infty, +\infty]$ будем сопоставлять некоторое ее <<{\it поднятие\/}>> 
$N^{\uparrow}\colon D\to [-\infty,+\infty]$, а именно: {\it для каждого $z\in D$ полагаем}
\begin{enumerate}
\item Если $\CC_{\infty}\setminus \clos D\neq \varnothing$ (непусто), то полагаем 
\begin{equation}\label{df:Nd}
N^{\uparrow}(z):= B(z,d(z);N)+\ln\frac{1}{d(z)}.
\end{equation}
\item Если $D=\CC$ --- комплексная плоскость, то для любого сколь угодно большого числа  $P>0$ можем положить
\begin{equation}\label{df:Nd_C}
N^{\uparrow}(z):= B\Bigl(z,\frac{1}{(1+|z|)^P};N\Bigr).
\end{equation}

\end{enumerate}	
\begin{remark} {\rm Обратим внимание, что предложенные функции-поднятия $N^{\uparrow}$, в отличие от предложенных ранее четырех таких функций (2a)--(2d) в \cite[перед теоремой 1]{Kha15}, которые здесь не приводятся ввиду громоздкости, значительно тоньше, компактнее и существенно более медленного роста.   Таким образом, предлагаемая ниже теорема 1 намного более точная, чем предшествующая ей \cite[теорема 1]{Kha15} 2015 г.}
\end{remark}

\begin{theorem}\label{th:gene} Пусть $D$ --- область в $\CC$, функции $N,M \in \sbh (D)$,
  $M-N\in \sbh (D)$ с мерой Рисса $\nu_{M-N}$, $N,M\neq \boldsymbol{-\infty}$, $\Lambda$ --- последовательность точек в $D$.
  
 Если $\Lambda$ --- последовательность неединственности для $\Hol(D;\exp N]$,  
то $n_\Lambda+\nu_{M-N}$ --- подмера для класса $\sbh(D;M]$.
 
Обратно, если $n_\Lambda+\nu_{M-N}$ --- подмера для класса $\sbh(D;M]$, $N$ --- непрерывная функция на $D$, то  последовательность точек $\Lambda$ --- последовательность неединственности для пространства $\Hol(D;\exp N^{\uparrow}]$ с подходящей весовой функцией-под\-н\-я\-т\-и\-ем $N^{\uparrow}$ из   
\eqref{df:Nd} при $\CC_{\infty}\setminus \clos D\neq \varnothing$ 
и с функцией-поднятием $N^{\uparrow}$ при произвольном фиксированным числе $P>0$ из   \eqref{df:Nd_C} при $D=\CC$.  
\end{theorem}

\begin{proof} {\rm Пусть $\Lambda$ --- последовательность неединственности для $\Hol(D;\exp N]$. Это означает,  что для функции $f_{\Lambda}$ с $\Zero_{f_\Lambda}=\Lambda$ найдется ненулевая функция $h\in \Hol(D)$, для которой произведение $f_{\Lambda}h\in 
\Hol (D;N]$, или $\log|f_\Lambda|+\log|h|\leq N$.	
  Введем обозначение
  \begin{equation}\label{sMN}
	s_{\nu_{M-N}}:=M-N\in \sbh (D)
\end{equation}
Тогда 
$	\log|f_\Lambda|+\log|h|+s_{\nu_{M-N}}\leq N+(M-N)=M \; \text{на $D$}$,
т.е. для меры Рисса $n_\Lambda+\nu_{M-N}$ субгармонической функции $\log |f_\Lambda|+s_{\nu_{M-N}}$ нашлась субгармоническая функция $v=\log|h|\neq \boldsymbol{-\infty}$, для которой сумма 
$\bigl(\log|f_\Lambda|+s_{\nu_{M-N}}\bigr)+v$ принадлежит классу $\sbh(D;M]$. Таким образом, установлено, что  $n_\Lambda+\nu_{M-N}$ --- подмера для  класса $\sbh(D;M]$.
Обратно, пусть $n_\Lambda+\nu_{M-N}$ --- подмера для класса $\sbh(D;M]$. Это значит, что
найдется функция $w\in \sbh (D)$ с которой в обозначении \eqref{sMN}
\begin{equation*}
	\log |f_{\Lambda}|+s_{\nu_{M-N}}+w\leq M+\const \quad \text{\it  на $D$}.
\end{equation*}
Иначе
\begin{equation*}
	\log |f_{\Lambda}|+M-N+w-\const \leq M \quad \text{\it  на $D$},
\end{equation*}   
т.е. 
\begin{equation}\label{in:fb}
\log|f_\Lambda|+v\leq N 	\quad \text{ на $D$}
\end{equation}
для функции $v=w-\const\in \sbh(D)$, $v\neq \boldsymbol{-\infty}$. 
Будет использовано  }
\begin{propos}[{\rm \cite[следствие  3]{KhaBai16}}]\label{pr:d} Пусть $D$ --- область в $\CC$,  удовлетворяющая условию $\CC_{\infty}\setminus \clos D \neq \varnothing$ с непрерывной функцией $d\colon D\to (0,1]$, для которой выполнено условие \eqref{ddis}.
Для  субгармонической  в $D$ функции $v\neq \boldsymbol{-\infty}$
 найдется  ненулевая функция $h\in \Hol (D)$, удовлетворяющая условию
\begin{equation}\label{in:ub}
	\log |h(z)|\leq B(z,d(z);v)+\ln\frac{1}{d(z)} \quad\text{для всех $z\in D$}.
\end{equation}
\end{propos}
{\rm Применим усреднения по шарам к обеим частям 
\begin{equation*}
B\bigl(z,d(z);\log|f_{\Lambda}|\bigl)+B\bigl(z,d(z);v\bigl)\leq B\bigl(z,d(z);N\bigl),
\end{equation*}
или, в силу субгармоничности $\log |f|$, 
\begin{equation*}
\log|f_{\Lambda}|+B\bigl(z,d(z);v\bigl)\leq B\bigl(z,d(z);N\bigl).
\end{equation*}
Применяя к последнему неравенству соотношение \eqref{in:ub} предложения \ref{pr:d}, получаем требуемый случай с поднятием $N^{\uparrow}$ из \eqref{df:Nd}, поскольку ненулевая голоморфная функция $f_\Lambda h$ с подпоследовательностьб нулей $\Lambda$ принадлежит уже классу 
$\Hol(D;\exp N^{\uparrow}]$. 

Доказательство для случай $D=\CC$ проводится совершенно аналогично через}
 \begin{propos}[{\rm \cite[следствие  2 с комментарием]{KhaBai16}}]\label{pr:d_C}.
Для любой  субгармонической  в $\CC$ функции $v\neq \boldsymbol{-\infty}$
 для любого сколь угодно большого числа $P>0$ найдется  ненулевая целая функция $h\in \Hol (D)$, удовлетворяющая условию
\begin{equation}\label{df:Nd_Cv}
\ln |h(z)|\leq  B\Bigl(z,\frac{1}{(1+|z|)^P};v\Bigr) \quad\text{для всех $z\in \CC$}.
\end{equation}
\end{propos}

\end{proof}

\section{Положительно $\rho$-однородные\\ субгармонические функции}\label{gd}

Мы вынуждены  повторить и напомнить некоторые сведения, собранные в \cite[п.~3.1]{Kha15}
Пусть $\rho \in (0,+\infty)$.  Обозначим через $\rho\text{-\!}\shg (\CC)\subset \sbh (\CC)$ {\it множество
субгармонических положительно однородных при показателе\/ $\rho$  функций\/} $H\neq \boldsymbol{-\infty}$, т.\,е. $H(tz)=t^{\rho} H(z)$ при всех $z\in \CC$, $t\geq 0$. Через $\rho\text{-\!}\trc (\RR)$ обозначаем множество {\it $2\pi$-периодических\/ $\rho$-тригонометрически выпуклых\footnote{Используют также термины {\it тригонометрически\/ $\rho$-выпуклая,\/} или {\it тригонометрически выпуклая при показателе (порядке)\/ $\rho$.}} функций}\/ $h\colon \RR\to \RR$, \cite{Levin56}, \cite{Levin96}, \cite{Ev}, \cite{Mg} : 
\begin{equation}
	h(\theta)\sin \rho(\theta_2-\theta_1)\leq 	h(\theta_1)\sin \rho(\theta_2-\theta)+
	h(\theta_2)\sin \rho(\theta-\theta_1), \quad \theta_1\leq\theta\leq \theta_2<\theta_1+\frac{\pi}{\rho}\,.
\end{equation}
 Известно, что
 
\begin{enumerate}
 	\item\label{r1}  
	{Отображение-расширение $\ext\colon h\longmapsto \bigl(H\colon re^{i\theta}\mapsto h(\theta)r^{\rho}, 
	\; r\geq 0,\; \theta \in \RR\bigr)$, функций $h\colon \RR\to \RR$ задает аддитивную положительно однородную сохраняющую точную верхнюю грань биекцию выпуклого конуса\/  $\rho\text{-\!}\trc (\RR)$ на выпуклый конус\/ $\rho\text{-\!}\shg (\CC)$, и функции из 
	$\rho\text{-\!}\shg (\CC)$ и\/  $\rho\text{-\!}\trc (\RR)$ непрерывны\/}
		\cite{Levin56}--\cite[\S~2.3, I--VI]{Ev}, \cite[Свойство 9.5, Теоремы 9.12]{Mg};
\item\label{r2}  {\it функция\/ $H\in \rho\text{-\!}\shg (\CC)$ удовлетворяет локальному условию Липшица в  форме\/}
{\large(}см. \cite[\S~2.3, IV]{Ev} и детальнее  \cite[Свойство 9.25 и Следствие 9.26 с доказательством]{Mg}{\large)}
\begin{equation}\label{co:L}
	\bigl| H(z)-H(w)\bigr|\leq \rho \max_{\varphi\in \RR} H(e^{i\varphi})  \cdot \bigl(\max\bigl\{|z|,|w|\bigr\}\bigr)^{\rho-1} 	|z-w|, \quad z,w\in \CC,
	\end{equation}
	{\it и,\/}  как следствие из ($\rho$\ref{r1}), {\it функция\/ $h\in \rho\text{-\!}\trc (\RR)$  удовлетворяет условию Липшица}
		\begin{equation}\label{co:l}
		\bigl|h(\theta)-h(\vartheta)\bigr|\leq \rho \max_{\varphi \in \RR } h(\varphi) \cdot |\theta-\vartheta|, \quad \theta,\vartheta \in \RR;
\end{equation}
	\item\label{r3} в обозначениях из ($\rho$\ref{r1}) {\it плотность меры Рисса\/ $d\nu_H$ функции\/ $H\in \rho\text{-\!}\shg (\CC)$ в полярных координатах определяется как произведение плотностей мер 
	\begin{equation}\label{df:densR}
d\nu_H(re^{i\theta})=r^{\rho-1}\,dr\otimes \frac{1}{2\pi}\,\bigl(h''(\theta)+\rho^2h(\theta)\bigr)\,d \theta, \quad re^{i\theta}\in \CC, \; r\geq 0, \; \theta \in \RR,
\end{equation}
где производные понимаются в смысле теории распределений, или обобщеных функций, а\/  $h''+\rho^2h\geq 0$  --- положительная $2\pi$-периодическая мера на $\RR$.}
\end{enumerate}

Пусть $h_1, h_2\in \rho\text{-\!}\trc (\RR)$. Воспользуемся  терминологией диссертации А.\,В.~Абанина \cite[\S~2.5]{Abd} 1995 г.,   широко используемой при исследовании абсолютно представляющих систем.  Называем функцию $h_1$ {\it $\rho$-выпукло дополнимой\/} до $h_2$, если  $h_2-h_1\in  \rho\text{-\!}\trc (\RR)$. В этом случае представляется естественным называть и  \large{(}см. ($\rho$\ref{r1})\large{)} функцию $H_1:=\ext h_1\in \rho\text{-\!}\shg(\CC)$ 
$\rho$-выпукло дополнимой до  $H_2:=\ext h_2\in \rho\text{-\!}\shg(\CC)$.
Из ($\rho$\ref{r3}) сразу следует, что $h_1$ {\it $\rho$-выпукло дополнима\/} до $h_2$, если и только если в смысле теории распределений, или обобщенных функций, $(h_2-h_1)''+\rho^2(h_2-h_1)\geq 0$, т.\,е. слева положительная $2\pi$-периодическая мера на $\RR$.  В частности, отсюда
\begin{enumerate}
	\item {\it если функция\/ $g\in \rho\text{-\!}\trc (\RR)$ дифференцируема на\/ $\RR$, т.\,е.\/ $g\in C^1(\RR)$, и 
		\begin{equation}\label{co:251}
	g'(\psi)-g'(\varphi)	\geq -c(\psi-\varphi), \quad 0\leq \varphi <\psi <2\pi, 
	\end{equation}
	где $0<c<\rho^2\min\bigl\{ g(\theta)\colon \theta \in [0,2\pi)\bigr\}$, то любая функция $h\in \rho\text{-\!}\trc (\RR)$ 
класса $C^1(\RR)$, для которой 
\begin{equation}\label{co:252}
h'(\psi)-h'(\varphi)	\leq C(\psi-\varphi), \quad 0\leq \varphi <\psi <2\pi,	\quad \text{$C$ --- постоянная},
\end{equation}
 $\rho$-выпукло дополнима до функции\/  $qg\in \rho\text{-\!}\trc (\RR)$ для любого числа\/}  (см. \cite[Лемма 2.5.1 и ее доказательство]{Abd}) 
\begin{equation}\label{df:q}
	q>\frac{C+\rho^2\max h}{\rho^2\min g -c};
\end{equation}

\item {\it если\/ $g\in \rho\text{-\!}\trc (\RR)$, 
	$\min \bigl\{ g(\theta)+g(\theta+\pi /\rho)\colon \theta \in [0, 2\pi]\bigr\}>0$,
и $g'$ --- неубывающая на\/ $[0,2\pi)$ {\rm \large(}в частности, сюда включается и случай постоянной функции\/ $g(\theta)\equiv R>0$, $\theta\in \RR${\rm \large)}, а  $h\in \rho\text{-\!}\trc (\RR)$ удовлетворяет\/ \eqref{co:252}, то функция\/ $h$ $\rho$-выпукло дополнима до\/ $qg$ для любого\/ $q$ из  \eqref{df:q} при  $c=0$} \cite[\S~2.5, Следствие 1]{Abd};
\item {\it если\/ $g\in \rho\text{-\!}\trc (\RR)$ удовлетворяет\/ \eqref{co:251}, а  $h\in \rho\text{-\!}\trc (\RR)$ имеет ограниченную вторую производную, т.\,е. в \eqref{co:252} постоянная $C=\sup |h''|$, то  функция\/ $h$ $\rho$-выпукла  дополнима до $qg\in \rho\text{-\!}\trc (\RR)$\/ для любого числа\/ $q$ из\/} \eqref{df:q} \cite[\S~2.5, Следствие 1]{Abd}.
\end{enumerate}

\begin{theorem}\label{th:2}  Пусть функция $h_1\in \rho\text{-\!}\trc (\RR)$ $\rho$-выпукло дополнима до $h_2\in \rho\text{-\!}\trc (\RR)$, т.\,е. в обозначениях из\/ {\rm ($\rho$1)} функция $H_1:=\ext h_1\in \rho\text{-\!}\shg(\CC)$ 
$\rho$-выпукло дополнима до  функции $H_2:=\ext h_2\in \rho\text{-\!}\shg(\CC)$, а мера $\nu$ определена через произведение плотностей мер в полярных координатах по правилу\/ {\rm (см. и ср. с \eqref{df:densR})}
\begin{equation}\label{df:densRH}
d\nu (re^{i\theta})=r^{\rho-1}\,dr\otimes \frac{1}{2\pi}\,
\bigl((h_2-h_1)''+\rho^2(h_2-h_1)\bigr) (\theta)\,d \theta, \quad re^{i\theta}\in \CC, \; r\geq 0, \; \theta \in \RR.
\end{equation}
Далее, пусть $\Lambda$ последовательность точек в\/ $\CC$. Тогда
если $\Lambda$ --- последовательность неединственности для\/ $\Hol(\CC;\exp H_1]$,  
то\/ $n_\Lambda+\nu_{M-N}$ --- подмера для класса\/ $\sbh(\CC;H_2]$, и обратно, если 
$n_\Lambda+\nu_{M-N}$ --- подмера для класса $\sbh(\CC;H_2]$, то
 при любом значении $\rho >0$ последовательность $\Lambda$  --- последовательность неединственности для $\Hol(\CC;\exp H_1]$.
\end{theorem}
Доказательство опускаем. Оно почти дословно повторяет доказательство \cite[теорема 2]{Kha15}, в которой различаются случаи $\rho \leq 1$ и $\rho >1$, где в заключении в последнем случае возникала существенная добавка-мультипликатор в виде   $\Hol(\CC; p\exp H_1]$, где $p$ достаточно быстро растущий многочлен. Такой результат со степенной добавкой значительно ослабляет теорему  \ref{th:2}.
Улучшение  в данном случае достигается за счет очень мало отличающейся функции-поднятия
из \eqref{df:Nd_C} от исходной функции $N$.
\begin{remark} {\rm Возможны обобщения результатов статьи на функции нескольких комплексных переменных, что предполагается проделать в ином месте.}
\end{remark}

\end{document}